\documentclass[11pt]{amsart}
\usepackage{amsmath}
\usepackage{amsfonts}
\usepackage{amssymb}
\usepackage{amsthm}
\usepackage[margin=.8in]{geometry}
\usepackage[dvipsnames]{xcolor}

\usepackage{tikz}
\usetikzlibrary{calc, positioning, fit, shapes.misc}
\usetikzlibrary{shapes.geometric}

\makeatletter
\def\@settitle{\begin{center}\LARGE\@title\end{center}}
\makeatother
\usepackage{float}

\newcommand{\GG}{\mathbb G}
\newcommand{\G}{\mathcal G}
\usepackage{titlesec}

\makeatletter
\def\@settitle{\begin{center}\Large\textbf{\@title}\end{center}}

 \renewenvironment{abstract}{
 \begin{center}%
 {\vspace{\baselineskip}
 \normalfont\fontsize{12}{12}\bfseries\abstractname
 \vspace{\baselineskip}}
 \end{center}\quotation
 }

 \titleformat{\section}
   {\normalfont\fontsize{14}{16}\bfseries}
   {\thesection}
   {1em}
   {}

 \titleformat{\subsection}
   {\normalfont\fontsize{12}{14}\bfseries}
   {\thesubsection}
   {1em}
   {}
  
 \titleformat{\subsubsection}
   {\normalfont\fontsize{12}{12}\bfseries}
   {\thesubsubsection}
   {1em}
   {}

\makeatother

\usepackage{enumitem}
\usepackage{bm}

\usepackage{tikz}
\usepackage{graphicx}

\usepackage[breaklinks=true]{hyperref} 

\usepackage{breakcites}
\hypersetup{colorlinks=true, citecolor=red, linkcolor=blue}

\renewcommand{\ge}{\geqslant}
\renewcommand{\le}{\leqslant}

\newcommand{\Ga}{\Gamma}

\newcommand{\ga}{\gamma}

\renewcommand{\L}{\mathcal L}
\newcommand{\Ind}{\mathcal I}

\newcommand{\CS}{\text{CSpan}}
\newcommand{\p}{\textbf p}

\newcommand{\E}[1]{\mathbb E\left[#1\right]}

\renewcommand{\P}[1]{\mathbb{P}\left(#1\right)}
\newcommand{\Exp}[1]{\exp\left(#1\right)}
\renewcommand{\O}[1]{O\left(#1\right)}
\renewcommand{\o}[1]{o\left(#1\right)}
\newcommand{\Log}[1]{\log\left(#1\right)}

\newcommand{\lam}{\lambda}

\newcommand{\C}{\mathcal C}

\newcommand{\F}{\mathcal F}
\renewcommand{\L}{\mathcal L}

\newcommand{\ep}{\varepsilon}

\newcommand{\R}{\mathbb R}

\renewcommand{\k}{\kappa}

\newcommand{\Mu}[1]{\mu\left(#1\right)}

\newcommand{\I}[1]{\mathbf 1_{\left\{#1\right\}}}

\newtheorem{tm}{Theorem}[section]

\newtheorem{lm}[tm]{Lemma}
\newtheorem{co}[tm]{Corollary}

\newtheorem{re}[tm]{Remark}

\theoremstyle{definition}
\newtheorem{df}[tm]{Definition}

\title{Asymptotic linearity of binomial random hypergraphs
via \\
cluster expansion under graph-dependence}
\date{}

\begin{document}

\maketitle

\begin{center}
Rui-Ray Zhang \\
rui.zhang@monash.edu \\
School of Mathematics, Monash University
\end{center}

\begin{abstract}

Let integer $n \ge 3$ and integer $r = r(n) \ge 3$.
Define the binomial random $r$-uniform hypergraph $H_r(n, p)$ 
to be the $r$-uniform graph on the vertex set $[n]$
such that each $r$-set is an edge independently with probability $p$.
A hypergraph is linear if every pair of hyperedges {intersects} in at most one vertex.
We study the probability of linearity of random hypergraphs $H_r(n, p)$
via cluster expansion {and give more precise asymptotics of the probability in question, 
improving the asymptotic probability {of linearity} obtained by McKay and Tian,
{in particular, when $r=3$ and $p = o(n^{-7/5})$}}.
\end{abstract}

\section{Introduction}

{
For all positive integer {$n \ge 3$}, let $[n]$ denote {the} integer set $\{ 1, 2, \ldots, n \}$.
Let $[n]_t = n (n - 1) \cdots (n - t + 1)$ denote the {$t$-th}
falling factorial for {every} positive integer $t \le n$.
}
Define the binomial random $r$-uniform hypergraph $H_r(n, p)$ 
to be the $r$-uniform hypergraph ($r$-graph for short) on {the} vertex set $[n]$
such that each $r$-element subset
($r$-set for short) is an edge independently with probability $p$.
A hypergraph is linear if every pair of hyperedges {intersects} in at most one vertex.
Let $\L_r(n)$ be the set of all linear $r$-uniform hypergraphs with $n$ vertices.
All asymptotics in this note are with respect to $n \to \infty$.
{We study the probability of linearity of random hypergraphs $H_r(n, p)$,
and improve the following result 
by giving more precise asymptotics of the probability.}

\begin{tm}[{\cite[Theorem 1.2]{mckay2020asymptotic}}]
Let $r = r(n) \ge 3$.
If $ p \binom{n}{r} = {\O{ r^{-2} n }}$, then
\begin{align*}
\P{ H_r(n, p) \in \L_r(n) }
= \Exp{ 
- \frac{ [r]_2^2 }{4 n^2 } \binom{n}{r}^2 p^2 
+ \O{ \frac{ r^6 }{ n^3 } \binom{n}{r}^2 p^2 } }.
\end{align*}
If {$r^{-2} n \le p \binom{n}{r} = \o{ r^{-3} n^{3/2} }$}, then
\begin{align*}
\P{ H_r(n, p) \in \L_r(n) }
= \Exp{ 
- \frac{ [r]_2^2 }{4 n^2 } \binom{n}{r}^2 p^2
+ \frac{ (3r-5) [r]_2^3 }{ 6n^4 } \binom{n}{r}^3 p^3
+ \O{ \frac{ \log^3 ( r^{-2} n ) }{ \sqrt{ \binom{n}{r} p } }
+ \frac{ r^6 }{ n^3 } \binom{n}{r}^2 p^2 } }.
\end{align*}
\end{tm}
For random $3$-uniform hypergraphs,
the above theorem gives that if $ p = \o{ n^{-3/2} }$, 
then
\begin{align}
\P{ H_3(n, p) \in \L_3(n) }
= \Exp{ - \frac{1}{4} n^4 p^2
+ \frac{2}{3} n^5 p^3 
+ \o{1}}.
\label{mc-3}
\end{align}

%
%

\subsection{Cluster expansion and dependency graphs}

{
Cluster expansion is a powerful tool in the rigorous study of statistical mechanics.
It was pioneered by Mayer in the 1930's and remains widely used nowadays,
see, for example, \cite{friedli2017statistical}.
The cluster expansion allows {us} to express the logarithm of partition function as a sum over clusters. 
Here we introduce the standard cluster expansion setting,
which is formulated in a way that is convenient for our application.}
\begin{df}
Given an undirected graph $G = (V(G), E(G))$.
\begin{enumerate}[leftmargin=1cm, label=(d\arabic*)]

\item A connected component of $G$ is a maximal set of vertices 
such that every pair of vertices is connected by a path.
The number of connected components of $G$ is denoted by $c(G)$.

\item {The set of polymers $\C(G)$ of $G$ is
the vertex sets of all connected induced subgraphs of $G$, namely,
\begin{align*}
\C(G) = \{ C \subseteq V(G): c(G[C]) = 1 \},
\end{align*}}
where $G[C]$ denotes the subgraph of $G$ induced by {the} vertex set $C$.
{For any two distinct polymers $C_i, C_j \in \C(G)$, 
we write $C_i \sim C_j$ if $C_i \cup C_j \in \C(G)$;
otherwise, $C_i \not\sim C_j$.
Equivalently, $C_i \not\sim C_j$ if $d_G(C_i, C_j) > 1$ and otherwise, $C_i \sim C_j$,
where $d_G(\cdot, \cdot)$ denotes the graph distance in $G$,
that is, the number of edges in the shortest path among two subsets of vertices.
Note that if we have that $C_i \sim C_j$
and $C_i \cap C_j = \emptyset$,
then $C_i$ and $C_j$ are adjacent in $G$,
that is, 
there exists an edge in $E(G)$, with one endpoint in $C_i$ and the other in $C_j$
}.
The size of a polymer, {denoted by $|C|$}, is the number of vertices in it.

\item 
{
For every non-empty ordered multiset of polymers $( C_1, \ldots, C_n ) \in \C(G)^n$, 
let $\GG( C_1, \ldots, C_n )$ be the graph on 
$[n]$ with $\{ i, j \} \in E( \GG )$ if $ C_i \sim C_j $.
}

{For instance, fix a polymer $C\in \C(G)$.
For a multiset of $n$ copies of $C$,
we have $\GG( C, \ldots, C ) = K_n$,
where $K_n$ denotes the complete graph on $[n]$}.

\item A cluster $\ga$ is a non-empty ordered multiset of polymers 
$( C_1, \ldots, C_{ |\ga| } )$ such that {$\GG(\ga) = \GG( C_1, \ldots, C_{ |\ga| } )$} is connected.
The size of a cluster $\ga$, {denoted by $|\ga|$}, is the number of polymers in it,
and the number of vertices of a cluster $\ga$, {denoted by $\|\ga\|$},
is the sum of size of polymers it contains, {that is, $\|\ga\| = \sum_{C \in \ga} |C| $}.

\item The set of all clusters of $G$ is denoted by $\Ga( G )$.
The set of all clusters of $G$ with pairwise disjoint polymers is 
denoted by
\begin{align*}
{ \Ga_{\emptyset}( G ) = 
\left\{ \ga \in \Ga( G ):
C_i \cap C_j = \emptyset \text{ for any distinct } C_i, C_j  \in \ga \right\}.}
\end{align*}
{Note that each element in $\Ga_{\emptyset}( G )$ is a cluster 
whose elements form a partition of a polymer,
since for every $\ga \in \Ga_{\emptyset}( G )$, polymers 
{$\{ C: C \in \ga \}$} are disjoint and their union $\cup_{C \in \ga} C \in \C(G)$.}

\item The set of all connected spanning subgraphs of $G$ is 
\begin{align*}
\CS(G) = \{ (V(G), E) : E \subseteq E(G), c((V(G), E)) = 1 \}.
\end{align*} 
{Here for every graph $H \in \CS(G)$,
we have $V(H) = V(G)$, $E(H) \subseteq E(G)$, and $c(H) = 1$}.
\end{enumerate} 
\end{df}

{The cluster expansion method
can be naturally combined with
dependency graphs.
The dependency graph models have been widely used in probability and statistics
to establish normal or Poisson approximation
via the Stein's method, cumulants, etc.
(see, for example, \cite{janson1988normal, janson1990poisson}).
They are also heavily used in probabilistic combinatorics, 
such as Lov\'asz local lemma \cite{erdos1975problems}, 
Janson's inequality \cite{janson1988exponential}, 
concentration inequalities \cite{zhang2022janson}, etc.}

Given a graph $G = (V, E)$,
we say that random variables $\{ X_i \}_{i \in V}$ are $G$-\textit{dependent}
if for any disjoint $S, T \subset V$ {such that $d_G(S, T) > 1$},
random variables $\{ X_i \}_{i \in S}$ and $\{ X_j \}_{j \in T}$ are independent.
{
Or equivalently, 
random variables $\{ X_i \}_{i \in C_1}$ and $\{ X_j \}_{j \in C_2}$ are independent
for any two distinct polymers $C_1$ and $C_2$ of the graph $G$ such that $C_1 \not\sim C_2$.}

Note that the dependency graph for a set of random variables 
may not be necessarily unique
and the sparser ones are the more interesting ones.
{Since no two disjoint $S, T \subset [n]$ are non-adjacent in $K_n$,
then the trivial dependency graph $K_n$ is a valid dependency graph 
for any set of variables $\{ X_i \}_{i \in [n]}$.}

Given $G$-dependent random variables $\{ X_i \}_{i \in V(G)}$,
for every set of vertices $S \subseteq V(G)$,
the \textit{joint moment} of random variables $\{ X_i \}_{i \in S}$ is defined by 
\begin{align*}
\Mu{ S } = \E{ \prod_{i \in S} X_i },
\end{align*} 
with $\mu(\emptyset) := 1$.
{
Let $\{ C_i \}_{i \in [n]}$ be a set of pairwise non-adjacent disjoint polymers of $G$,
in other words, for all distinct $i, j \in [n]$, we have $C_i \not\sim C_j$, 
or equivalently, $d_G(C_i, C_j) > 1$.
Then one important factorisation property 
for $G$-dependent variables, 
following from the definition of dependency graph, is that
\begin{align}
\Mu{ \bigcup_{i \in [n]} C_i } = \prod_{i \in [n]} \Mu{ C_i }.
\label{fact}
\end{align}
}

{
Let $\{ X_v \}_{v \in V(G)}$ be $G$-dependent random indicators and $X = \sum_{v \in V(G)} X_v$.
In our application, 
each indicator indicates the occurrence of some combinatorial structure,
and indicators are dependent with a dependency graph.
By writing the probability of the non-existence of some combinatorial structure
$\P{X = 0}$ as a partition function, 
the cluster expansion then
gives the formal expansion formula as a sum over clusters, 
whose truncation approximates the asymptotic probability.
This is inspired by \cite{scott2005repulsive},
in which they also treat $\P{X = 0}$ as a partition function and investigate the connections
between cluster expansion and the Lov{\'a}sz local lemma,
giving {a} lower bound for $\P{ X = 0 }$.}

The standard cluster expansion gives the formal cluster expansion
\begin{align}
\log \P{ X = 0 }
= \sum_{ \ga \in \Ga( G ) }
\frac{ \phi(\ga) }{ |\ga|! } (-1)^{\|\ga\|}
\prod_{C \in \ga} \Mu{C},
\label{ce}
\end{align}
with \textit{Ursell function} $\phi: \Ga( G ) \rightarrow \R$ defined by
\begin{align}
\phi(\ga)
= \sum_{ H \in \CS( \GG( \ga ) ) } (-1)^{e_H},
\label{ursell}
\end{align}
where $e_{H}$ denotes the number of edges of {the} graph $H$.
Note that 
if the cluster $\ga$ contains one single polymer $C \in \C(G)$, 
then $\phi(\ga) = 1$,
{because} $\GG( C ) = K_1$.


{For completeness, we include a simple derivation of {equation} \eqref{ce},
following the routine cluster expansion derivation procedure
(see, for example, \cite[Section 2.2]{scott2005repulsive}
or \cite[Proposition 5.3.]{friedli2017statistical}).}
First, the inclusion-exclusion formula gives
\begin{align}
\P{ X = 0 }
= \sum_{S \subseteq V(G)} (-1)^{|S|} \Mu{ S }.
\label{in-ex}
\end{align}

{Let $\G_{c}$ be a graph on vertex set $\C(G)$
such that for all distinct $C_i, C_j \in \C(G)$,
if $ C_i \sim C_j $,
then $\{ C_{i}, C_{j} \} \in E(\G_{c})$.
Next we utilise the factorisation property
{as shown in \eqref{fact} to
prove that the right hand side of equation (\ref{in-ex})} can be written as 
some partition function of hard-core model,
more specifically,
as a summation over independent sets of graph $\G_{c}$,}
\begin{align}
\P{ X = 0 } 
= \sum_{U \in \Ind(\G_{c}) } \prod_{C \in U} (-1)^{|C|} \Mu{C},
\label{part}
\end{align}
where $\Ind(G)$ denotes the set of all independent sets for every graph $G$.

For every $S \subseteq V(G)$ such that $S \in \C(G)$, we have {$ \{ S \} \in \Ind(\G_{c})$}.
For every $S \subseteq V(G)$ such that $S \notin \C(G)$, 
we have $S$ {induces} a union of pairwise non-adjacent maximal connected subgraphs,
that is, 
there exists a unique set of polymers
$U \in \Ind(\G_{c})$
such that $S = \cup_{C \in U} C$,
and $C_i \not\sim C_j$ for all pairs of distinct $C_i, C_j \in U$.
{The factorisation property \eqref{fact}} then gives
\begin{align}
(-1)^{|S|} \Mu{ S } = \prod_{C \in S} (-1)^{|C|} \Mu{ C }.
\label{decom}
\end{align}
Conversely, for every $U \in \Ind(\G_{c})$, we have $\cup_{C \in U} C \subseteq V(G)$,
thus $U$ determines $S$ uniquely;
combining with {equation \eqref{decom}, 
it follows that equations (\ref{in-ex}) and (\ref{part}) are equivalent}.

Now we derive the formal cluster expansion. 
{Let $\binom{S}{i}$ denote the family of {$i$-sets} of $S$
for every set $S$ and every positive integer $i \le |S|$.}
From {equation} \eqref{part}, we have
\begin{align*}
\P{ X = 0 } 
&= \sum_{U \subseteq \C(G)} \prod_{C \in U} (-1)^{|C|} \Mu{ C } 
\prod_{ \{ C_i, C_j \} \in \binom{U}{2} } \I{ C_i \not\sim C_j } \\
&= 1 + \sum_{n \ge 1} \frac{1}{n!} 
\sum_{ (C_1, \ldots, C_n) \in \C(G)^{n} }
\prod_{i \in [n]} (-1)^{|C_i|} \Mu{ C_i } 
\prod_{ 1 \le i < j \le n } \I{ C_i \not\sim C_j }.
\end{align*}
Note that by a simple expansion,
\begin{align*}
\prod_{ 1 \le i < j \le n } \I{ C_i \not\sim C_j }
= \prod_{ 1 \le i < j \le n } \left( 1 - \I{ C_i \sim C_j } \right)
= \sum_{ H \in \mathfrak G_n } (-1)^{e_H} \prod_{ \{ i,j \} \in E(H) } \I{ C_i \sim C_j },
\end{align*}
where $\mathfrak G_n$ denotes the set of all graphs on $n$ vertices. Then formally, we have
\begin{align*}
\P{ X = 0 } 
= 1 + \sum_{n \ge 1} \frac{1}{n!} 
\sum_{ H \in \mathfrak G_n } W(H),
\end{align*}
{where
\begin{align*}
W(H)
= \sum_{ (C_1, \ldots, C_{v_H}) \in \C(G)^{v_H} }
(-1)^{e_H} \prod_{ \{ i,j \} \in E(H) } \I{ C_i \sim C_j }
\prod_{k \in [v_H]} (-1)^{|C_k|} \Mu{ C_k },
\end{align*}
{and $W(H)$} satisfies
\begin{enumerate}[leftmargin=1cm, label=(a\arabic*)]
\item $W(H) = W(H')$ whenever $H$ and $H'$ are isomorphic $H \cong H'$, 
that is, differ only by vertices relabelling;
\item $W(H) = W(H_1) W(H_2)$ whenever $H$ is isomorphic to the disjoint
union of $H_1$ and $H_2$.
\end{enumerate}}
Let $\mathfrak C_n$ be the set of all connected graphs on $n$ vertices.
{Via the exponential formula {\cite[Corollary 5.1.6]{stanley1999enumerative}}}, we 
reduce the sum over the set of all graphs to the set of all connected graphs
\begin{align*}
\log \P{ X = 0 }
= \sum_{n \ge 1} \frac{1}{n!} 
\sum_{ H \in \mathfrak C_n } W(H) 
= \sum_{ \ga \in \Ga( G ) }
\frac{ 1 }{ |\ga|! } \sum_{ H \in \CS( \GG( \ga ) ) } (-1)^{e_H + \|\ga\|}
\prod_{C \in \ga} \Mu{C},
\end{align*}
where $\Ga( G )$ denotes the set of all clusters of $G$.
{A similar derivation of the cluster expansion 
utilizing the exponential formula also appears in \cite[Proposition 5.3]{friedli2017statistical}}.
Then {equation} \eqref{ce} follows.

\begin{re}

\begin{itemize}
\item[(r1)]
{In probability theory and statistical physics,}
given a graph $H$ and a vector $\p = \{ p_v \}_{v \in V(H)}$,
the partition function of the hard-core model (also the independence polynomial) 
on $H$ is defined by
\begin{align*}
{ \mathbf I_H(\p) } = \sum_{U \in \Ind(H)} \prod_{i \in U} p_{i}.
\end{align*}
{The cluster expansion is essentially the multivariate Taylor series for 
{$\log \mathbf I_H(\p)$ in variables $\{ p_v \}_{v \in V(H)}$} around $\bm 0$.}
Let $\bm\mu = ( (-1)^{|C|} \mu(C))_{C\in \C(G)}$. Then
{equation (\ref{part})} can be regarded as the partition function $\mathbf I_{\G_c}(\bm\mu )$ 
of the hard-core model on $\G_c$.

\item[(r2)]
{For independent indicators $\{ X_i \}_{i \in [n]}$,
if $0 \le \E{ X_i } < 1$ for all $i \in [n]$,
then we have the Taylor series of logarithmic function
\begin{align}
\log \P{ X = 0 }
= \sum_{i \in [n]} \Log{ 1 - \E{ X_i } }
= - \sum_{i \in [n]} \sum_{j \ge 1} \frac{1}{j} \E{ X_i }^j.
\label{taylor}
\end{align}

The empty graph $\overline K_{n} := ([n], \emptyset)$ is a valid dependency graph 
for this independent case.
Since the polymers of  $\overline K_{n}$ are all of size one containing a single vertex,
and the clusters of $\overline K_{n}$ are all multisets containing 
multiple copies of the same vertex,
then for independent indicators $\{ X_v \}_{v \in [n]}$, expansion {in equation (\ref{ce})} becomes}
\begin{align}
\log \P{ X = 0 }
= \sum_{ \ga \in \Ga[ \overline K_{n} ] }
\frac{ \phi(\ga) }{ |\ga|! }
(-1)^{\|\ga\|} \prod_{C \in \ga} \Mu{C}
= \sum_{i \in [n]} \sum_{j \ge 1} 
\frac{1 }{j!} \sum_{ H \in \CS(K_j) } (-1)^{e_H} (-1)^j \E{ X_i }^j.
\label{in-exp}
\end{align}
Comparing {\eqref{taylor}
and \eqref{in-exp}, 
it follows that}
\begin{align}
\sum_{ H \in \CS(K_n) } (-1)^{e_H} 
= (-1)^{n-1} (n-1)!,
\label{cs-id-k}
\end{align}
which is well-known, see, 
for example, {\cite[Eq. (2.13)]{scott2005repulsive}} or {\cite[Eq. (3.37)]{stanley2011enumerative}}.
\end{itemize}
\end{re}

\section{Linearity of binomial random hypergraphs}


{
The probability of random hypergraphs being linear
is equal to the probability of the non-existence of hyperedge pairs intersecting in more than one vertex.
Then we accordingly define a set $\F$ of `forbidden' hypergraphs,
containing all $r$-graphs $( e_1 \cup e_2, \{ e_1, e_2 \} )$
on vertex set $e_1 \cup e_2$
such that 
$2 \le | e_1 \cap e_2 | < r$:
\begin{align}
\F = 
\bigcup_{2 \le t \le r - 1} 
{\Big\{} ( e_1 \cup e_2, \{ e_1, e_2 \} ):
|e_1| = |e_2| = r, | e_1 \cap e_2 | = t {\Big\}}.
\label{f}
\end{align}
Note that removing isomorphic duplicates from $\F$ does not 
affect the probability that we are interested in, 
we thus assume that the $r$-graphs in $\F$ are pairwise non-isomorphic. 
We also can assume that {no hypergraphs in $\F$ have} isolated vertices.

{The complete $r$-graph on $n$ vertices, denoted by $K^r_n$, is 
the hypergraph consisting of n vertices and all possible edges of size $r$, 
that is, $K^r_n = ([n], \{S \subseteq [n] : |S| = r\})$.}
For every $F \in \F$,
let $A^F$ be the set of all subgraphs 
of $K^r_n$ that are isomorphic to $F$.
There are $[n]_{v_F} / \text{Aut}(F)$ such subgraphs, 
where $\text{Aut}(F)$ denotes the number of automorphisms of $F$.
Let $A^{\F} = \cup_{F \in \F} A^F$.
Then {the} random variable
\begin{align*}
X = \sum_{F \in A^\F} 
{\I{ F \subset H_r(n, p) }}
\end{align*}
counts the number of all copies of all forbidden $r$-graphs of $\F$ in $H_r(n, p)$.
Hence the probability of random hypergraphs being linear $\P{ H_r(n, p) \in \L_r(n) }$
equals the probability $\P{X = 0}$ {such that} $H_r(n, p)$ 
avoids all copies of all $r$-graphs of $\F$.

Next we define a dependency graph $D$ {on the vertex set $A^{\F}$ for random indicators 
{$ \{ \I{ F \subset H_r(n, p) } \}_{F \in A^{\F}}$}}
with two indicators being dependent if the corresponding forbidden $r$-graphs share edges,
\begin{align}
D = \left( A^{\F}, 
\left\{ \{ F_1, F_2 \} \in \binom{A^{\F}}{2}:
E(F_1) \cap E(F_2) \ne \emptyset \right\} \right).
\label{Dind}
\end{align}

Using the above dependency graph for random indicators of the forbidden structures,
we then can utilize the truncated cluster expansion series 
to approximate the probability of a binomial random $r$-uniform hypergraph being linear
and to obtain more precise asymptotic probability {of linearity}.
In this setting, a polymer is a set of forbidden subgraphs 
whose induced subgraph in $D$ is connected.
We will use the truncation of cluster expansion series 
involving clusters with certain restricted number of forbidden subgraphs.
For every integer $k>0$, denote the $k$-th term of the cluster expansion 
and the $k$-th truncated expansion with \textit{disjoint} polymers as
\begin{align*}
L^{\emptyset}_{D, k} := 
\sum_{ \ga \in \Ga_{\emptyset}( D ): \|\ga\| = k }
\frac{ \phi(\ga) }{ |\ga|! } (-1)^{\|\ga\|}
\prod_{C \in \ga} \Mu{C}
\quad\text{ and }\quad
T^{\emptyset}_{D, k} := \sum_{i \in [k-1]} L^{\emptyset}_{D, i}.
\end{align*}
}

Now we are ready to state our main result utilizing the truncated cluster expansion.
\begin{tm}
{Let $r = r(n) \ge 3$.}
If $p = \o{ n^{2 - r} }$,
then for every integer $k > 0$,
\begin{align}
\P{ H_r(n, p) \in \L_r(n) }
= \Exp{ T^{\emptyset}_{D, k} + \O{ \Delta_{k+1}(D) } + \o{1} },
\label{exp-del}
\end{align}
where $D$ is the dependency graph for the indicators of forbidden $r$-graphs defined by 
\eqref{Dind} and $\Delta_i(D)$ denotes 
the sum of joint moments over polymers of size $i$ in the graph $D$,
\begin{align*}
\Delta_i(D) = \sum_{C \in \C(D): |C| = i} \mu( {C} ).
\end{align*}
Moreover, for any $\ep > 0$,
if $p = \o{ n^{2 - r -\ep} }$,
then there exists an integer $k = k( \ep ) > 0$ such that
\begin{align}
\P{ H_r(n, p) \in \L_r(n) }
= \Exp{ T^{\emptyset}_{D, k} + \o{1} }.
\label{lprop}
\end{align}
\label{tm-linear}
\end{tm}

{Theorem \ref{tm-linear} gives the more precise asymptotics of the probability
of random hypergraphs being linear.
We next consider a specific example, 
by restricting to the $3$-uniform hypergraphs case,
and computing only the first few terms of the series explicitly
for illustration purpose.
This extends the asymptotic probability {of linearity for $H_3(n, p)$} 
given by McKay and Tian in \eqref{mc-3}.}

\begin{co}
If $p = { \o{ n^{-7/5} } }$, then
\begin{align}
\P{ H_3(n, p) \in \L_3(n) }
= \Exp{ - \frac{1}{4} n^4 p^2
+ \frac{2}{3} n^5 p^3 
- { \frac{55}{24} } n^6 p^4
+ \frac{3}{2} n^3 p^2
+ \o{1} }.
\label{linear}
\end{align}
\label{co-3linear}
\end{co}

\section{Proofs {of main results}}

The density of a graph $G$ is defined by $ d(G) = e_G / v_G $,
where $v_G$ and $e_G$ are the numbers of vertices and edges of $G$ respectively.
Another commonly used (see, for example, \cite{rucinski1988small, mousset2020probability}) 
density measure $m_\star(G)$ is defined by
\begin{align}
m_\star(G) = \min_{ H \subseteq G, e_H \ge 1 } \frac{ e_G - e_H }{ v_G - v_H }.
\label{ex den}
\end{align}

{Next we introduce joint cumulant, which is a fundamental tool in probability theory.}
Given $G$-dependent random variables $\{ X_v \}_{v \in V(G)}$,
for every set of vertices $S \subseteq V(G)$,
the \textit{joint cumulant} of random variables $\{ X_i \}_{i \in S}$ is defined by 
\begin{align}
\k( S )
= \sum_{\pi \in \Pi(S)} (-1)^{|\pi| - 1} ( |\pi| - 1 )! \prod_{P \in \pi} \Mu{ P },
\label{mo cu}
\end{align}
where $\Pi(S)$ denotes the set of all partitions of $S$.
{
The joint cumulant $\k( S )$ 
can be regarded as a measure of the {mutual dependences} of the variables in $S$. 
An important property of the joint cumulant $\k( S )$ is that 
if $S$ can be partitioned into
two subsets $S_1$ and $S_2$ such that 
{the} variables in $S_1$ are independent of {the} variables in $S_2$, then
$\k(S) = 0$.
In other words, if $S\not\in \C(G)$, then $\k( S ) = 0$
(see, for example, \cite{MR725217}).}

{ Given a family $\F$ of $r$-graphs, 
we consider the probability that 
$H_r(n, p)$ is $\F$-free, that is, 
it simultaneously avoids all copies of all $r$-graphs in $\F$. }
Let $i > 0$ be an integer,
define $\k_i(D)$ to be 
the sum of joint cumulants over polymers of size $i$ {in the} dependency graph $D$, namely,
\begin{align*}
\k_i(D) = \sum_{C \in \C(D): |C| = i} \k( {C} ).
\end{align*}
\begin{lm}[{\cite[Corollary 12]{mousset2020probability}}]
Let $\F$ be a finite family of $r$-graphs 
and $p = p(n) \in (0, 1)$ satisfy 
\begin{align}
n p^{m_\star(\F)} = \o{1}
\qquad \text{and}
\qquad 
n p^{2d(\F)} = \o{1},
\label{ex den}
\end{align}
{where
\begin{align*}
m_\star(\F) = \min_{G \in \F} m_\star(G)
\quad\text{ and }\quad
d(\F) = \min_{G \in \F} d(G).
\end{align*}}
Then, for every integer $k > 0$, we have
\begin{align}
\left| \log \P{ H_r(n, p) \text{ is } \F\text{-free} } 
- \sum_{i \in [k]} (-1)^i \k_i(D) \right| = \O{ \Delta_{k+1}(D) } + \o{1}.
\label{cumu_exp}
\end{align}
Moreover, if $n p^{m_\star(\F)} = n^{-\ep}$ for {some} $\ep > 0$,
then there exists an integer $k = k(\ep, \F) > 0$ such that $\Delta_{k+1}(D) = \o{1}$.
\end{lm}

Here we introduce another lemma relating the cluster expansion and cumulant series in
{equation} (\ref{cumu_exp}).

\begin{lm}[Cumulant-cluster lemma]
Let $\{ X_v \}_{v \in V(G)}$ be $G$-dependent random indicators
and $k > 0$ be an integer.
Then
\begin{align}
T^{\emptyset}_{G, k+1}
= \sum_{ C \in \C_k( G )} (-1)^{|C|} \k(C),
\end{align}
where 
{$\C_k( G ) = \{ C \in \C(G): |C| \in [k] \}$ denotes
the set of polymers with size at most $k$}.

\label{lm-cum-clu}
\end{lm}

Now we are ready to prove the main result.
\begin{proof}[Proof of Theorem \ref{tm-linear}]
Let $\F$ defined by {equation} \eqref{f} be the set of forbidden $r$-graphs.
Since each $r$-set of $[n]$ is an edge independently with probability $p$ in $H_r(n, p)$, 
{then for distinct subgraphs $F_1, F_2 \in A^\F$, 
indicators {$ \I{ F_1 \subset H_r(n, p) } $ and $\I{ F_2 \subset H_r(n, p) }$}}
are dependent if $E(F_1) \cap E(F_2) \ne \emptyset$.
Hence graph $D$ defined by {equation} \eqref{Dind} is 
a dependency graph for random indicators $\I{ F \subset H_r(n, p) }$ for $F \in A^\F$.

Next we verify conditions in {equation} \eqref{ex den}.
Since 
\begin{align*}
m_\star(\F)
= \min_{G \in \F} \min_{ H \subseteq G, e_H \ge 1 } \frac{ e_G - e_H }{ v_G - v_H }
= \frac{ 1 }{ \displaystyle\max_{G \in \F} \max_{ H \subseteq G, e_H \ge 1 } ( v_G - v_H ) }
= \frac{1}{r - 2}
\end{align*}
and
\begin{align*}
d(\F) = \min_{G \in \F} d(G)
= \min_{G \in \F} \frac{e_G}{v_G}
= \frac{2}{ \displaystyle\max_{G \in \F} v_G}
= \frac{1}{r - 1},
\end{align*}
then
$2d(\F) \ge m_\star(\F)$ for all $r \ge 3$.
Thus if $np^{1/(r-2)} = \o{1}$, then
\begin{align}
\left| \log\P{ H_r(n, p) \in \L_r(n) }
- T^{\emptyset}_{D, k} \right| = \O{ \Delta_{k+1} } + \o{1}.
\end{align}
Moreover, if $p = \o{ n^{-(r-2)-\ep} }$ for some $\ep > 0$, then
combining with Lemma \ref{lm-cum-clu}, we complete the proof. 
\end{proof}

What remains is to show Lemma \ref{lm-cum-clu}.
We introduce an auxiliary lemma for its proof.
\begin{lm}
For all connected graph $H$, we have
\begin{align}
\sum_{\pi \in \Pi( V(H) )}
\sum_{ G \in \CS( K_{|\pi|} ) } (-1)^{e_G}
\prod_{P \in \pi} \I{ P \in \Ind( H ) }
= \sum_{ G \in \CS( H ) } (-1)^{e_G}.
\label{ind-part-ursell}
\end{align}
\label{lm-ind-part-ursell}
\end{lm}

We first introduce the chromatic polynomial.
Given a graph $H$ and a positive integer $\lam$,
a (proper) $\lam$-colouring of $H$ is a map
$\Phi: V(H) \rightarrow [\lam]$ such that $\Phi(u) \ne \Phi(v)$ for all $\{ u, v \} \in E(H)$.
The \textit{chromatic polynomial} ${ P_H(\lam) }$ of $H$ is the number of $\lam$-colourings of $H$.

Given a graph $H$ and a positive integer $k$, 
a partition 
containing $k$ {subsets} $\{ V_1, \ldots, V_k \}$ of $V(H)$ is called
a $k$-independent partition of $H$ if for every $i \in [k]$, 
we have $V_i \ne \emptyset$ and $V_i \in \Ind(H)$.
Let $\alpha(H, k)$ count the $k$-independent partition of $H$.
Then we have {the chromatic polynomial in factorial form}
\begin{align}
{ P_H(\lam) } = \sum_{k=1}^{{v_H}} \alpha(H, k) [\lam]_k,
\label{ch1}
\end{align}
(see, for example, {\cite[Theorem 1.4.1]{dong2005chromatic}}).
An equivalent formula for ${ P_H(\lam) }$
{written as a polynomial in $\lam$,
known as the Whitney-Tutte-Fortuin-Kasteleyn representation}
(see, for example, {\cite[Eq. (A.11)]{lovasz2012large}} or \cite[Eq. (1.2)]{sokal2001bounds}) 
is
{
\begin{align}
{ P_H(\lam) } = \sum_{ E \subseteq E(H) } (-1)^{|E|} \lam^{c(E)},
\label{ch2}
\end{align}
where $c(E) = c(V(H), E)$ counts the number of the connected components of subgraph $(V(H), E)$
for every edge set $E \subseteq E(H)$.}

\begin{proof}[Proof of Lemma \ref{lm-ind-part-ursell}]
{By inspecting {equations \eqref{ch2} and \eqref{ursell}},
one observes that the Ursell function is the linear term of the chromatic polynomial
(this is also a well-known fact, see, for exmaple, \cite{abdesselam2009clustering}).}
{Then we have
the right hand side of {equation} \eqref{ind-part-ursell}}
\begin{align}
\sum_{ G \in \CS( H ) } (-1)^{e_G}
= \frac{\mathrm{d} { P_H(\lam) } }{ \mathrm{d} \lam } \Big|_{\lam=0}
= \frac{\mathrm{d}}{ \mathrm{d} \lam } \left( \sum_{k=1}^{v_H} \alpha(H, k) [\lam]_k \right) \Big|_{\lam=0}
= \sum_{k=1}^{v_H} \alpha(H, k) (-1)^{k-1} (k - 1)!.
\label{ursell-chrom}
\end{align}
{Using the the combinatorial identity obtained before in (\ref{cs-id-k}),
the left hand side of {equation} \eqref{ind-part-ursell} can be rewritten as}
\begin{align}
\sum_{\pi \in \Pi( V(H) )}
\sum_{ G \in \CS( K_{|\pi|} ) } (-1)^{e_G}
\prod_{P \in \pi} \I{ P \in \Ind( H ) }
= \sum_{\pi \in \Pi( V(H) )} (-1)^{ |\pi| - 1 } (|\pi| - 1)! \prod_{P \in \pi} \I{ P \in \Ind( H ) }.
\label{ind-part}
\end{align}
Notice that {the right hand side of equation} (\ref{ind-part}) is a sum of $|\pi|$-independent partition for any $\pi \in \Pi( V(H) )$.
Thus we have
\begin{align}
\sum_{\pi \in \Pi( V(H) )} (-1)^{ |\pi| - 1 } (|\pi| - 1)! \prod_{P \in \pi} \I{ P \in \Ind( H ) }
&= \sum_{k=1}^{v_H} \sum_{\pi \in \Pi( V(H) ): |\pi|=k} 
(-1)^{ k - 1 } (k - 1)! \prod_{P \in \pi} \I{ P \in \Ind( H ) } \nonumber \\
&= \sum_{k=1}^{v_H} \alpha(H, k) (-1)^{ k - 1 } (k - 1)!.
\label{ind-part-chrom}
\end{align}
Then combining {\eqref{ursell-chrom} and \eqref{ind-part-chrom}, we complete the proof}.
\end{proof}

\begin{proof}[Proof of Lemma \ref{lm-cum-clu}]

From the cluster expansion, we have
\begin{align}
T^{\emptyset}_{G, k+1} 
= \sum_{i \in [k]} L^{\emptyset}_{G, i}
&= \sum_{ \ga \in \Ga_{\emptyset}( G ): \|\ga\| \in [k] }
\frac{ \phi(\ga) }{ |\ga|! } 
\prod_{C \in \ga} (-1)^{|C|} \Mu{C} \nonumber \\ 
&= \sum_{ \substack{ (C_1, \ldots, C_n) \in \Ga_{\emptyset}( G ) \\  
\sum_{i\in[n]} |C_i| \in [k] } }
\frac{1}{n!} \sum_{ H \in \CS( \GG( C_1, \ldots, C_n ) ) } (-1)^{e_H}
\prod_{i \in [n]} (-1)^{|C_i|} \Mu{C_i} \nonumber \\
&= \sum_{ \substack{ \{ C_1, \ldots, C_n \} \in \Ga_{\emptyset}( G ) \\ 
\sum_{i\in[n]} |C_i| \in [k] } }
\sum_{ H \in \CS( \GG( C_1, \ldots, C_n ) ) } (-1)^{e_H}
\prod_{i \in [n]} (-1)^{|C_i|} \Mu{C_i},
\label{rhs}
\end{align}
{where the first summation in the last line is an abuse of notation,
and denotes the summation over (unordered) sets of polymers.}
From the definition of joint cumulants, we get
\begin{align*}
\sum_{ C \in \C_k( G )} (-1)^{|C|} \k(C)
&= \sum_{ C \in \C_k( G )} (-1)^{|C|} 
\sum_{\pi \in \Pi(C)} (-1)^{|\pi| - 1} ( |\pi| - 1 )! \prod_{P \in \pi} \Mu{ P } \\
&= \sum_{ C \in \C_k( G )}
\sum_{\pi \in \Pi(C)} (-1)^{|\pi| - 1} ( |\pi| - 1 )! \prod_{P \in \pi} (-1)^{|P|} \Mu{ P }.
\end{align*}
Combining with identity {in (\ref{cs-id-k}), it follows that}
\begin{align*}
\sum_{ C \in \C_k( G )} (-1)^{|C|} \k(C)
= \sum_{ C \in \C_k( G )}
\sum_{\pi \in \Pi(C)} \sum_{ H \in \CS( K_{|\pi|} ) } (-1)^{e_H}
\prod_{P \in \pi} (-1)^{|P|} \Mu{ P }.
\end{align*}

Fix an arbitrary polymer $C \in \C_k( G )$,
for every partition $\pi = \{ P_1, \ldots, P_{m} \} \in \Pi(C)$, 
by factorising into pairwise non-adjacent maximal connected subgraphs,
{with each induced by $C_1, \ldots, C_n$ respectively},
there exists a unique finest partition {consisting of only polymers}
$\pi' = \pi'(\pi) = \{ C_1, \ldots, C_n \} \in \Pi(C)$ 
such that 
\begin{enumerate}[leftmargin=1cm, label=(p\arabic*)]
\item $n \ge m$,
\item for all $i\in[n]$, we have $C_i \in \C(G)$,
\item $(C_1, \ldots, C_n) \in \Ga_{\emptyset}( G )$, and
\item $\prod_{P \in \pi} \Mu{ P } = \prod_{C \in \pi'} \Mu{ C }$.
\end{enumerate}

\begin{figure}[H]
\begin{tikzpicture}[scale=.8]

\node[circle,draw=black] (a1) at (-.5,4) {$1$};  
\node[circle,draw=black] (a2) at (1.3,3)  {$2$}; 
\node[circle,draw=black] (a3) at (0.2,5)  {$3$};  
\node[circle,draw=black] (a7) at (3,2.5) {$5$};  
\node[circle,draw=black] (a4) at (4,5)  {$4$};  
\node[circle,draw=black] (a6) at (5.5,3)  {$6$};  
\node[circle,draw=black] (a5) at (6.5,4.5)  {$7$};  

\draw (a1) -- (a2);  
\draw (a2) -- (a3); 
\draw (a5) -- (a6);  
\draw (a2) -- (a4);  
\draw (a4) -- (a6);  
\draw (a6) -- (a7);  
\draw (a2) -- (a7);  

\node[line width=0.5mm, dotted, ellipse, draw=black, rotate=-40, scale=.8] [fit=(a1) (a3)] {};
\node[line width=0.5mm, dotted, ellipse, draw=black, rotate=45, scale=.8] [fit=(a2) ] {};
\node[line width=0.5mm, dotted, ellipse, draw=black, rotate=-10, scale=.72] [fit=(a4) (a5)] {};
\node[line width=0.5mm, dotted, ellipse, draw=black, rotate=10, scale=.72] [fit=(a7) (a6)] {};
\end{tikzpicture}
\caption{ A polymer of size seven with a partition 
$ \{ \{1, 3\}, \{ 2 \}, \{4, 7\}, \{5, 6\} \} $
and the corresponding polymer partition 
$ \{ \{1 \}, \{ 3\}, \{ 2 \}, \{4 \}, \{ 7\}, \{5, 6\} \} $
such that 
$\Mu{1, 3}\Mu{2}\Mu{4, 7}\Mu{5, 6}
= \Mu{1}\Mu{2}\Mu{3}\Mu{4}\Mu{5, 6}\Mu{7}$.}
\end{figure}
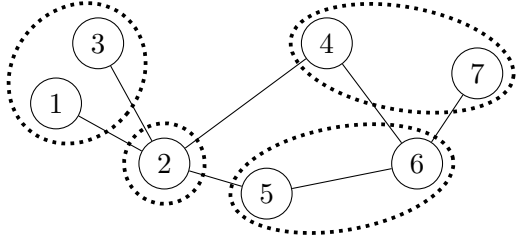
Then, we have
\begin{align*}
\sum_{ C \in \C_k( G )} (-1)^{|C|} \k(C)
&= \sum_{ C \in \C_k( G )}
\sum_{\pi \in \Pi(C)}
\sum_{ H \in \CS( K_{|\pi|} ) } (-1)^{e_H}
\prod_{C \in \pi'(\pi)} (-1)^{|C|} \Mu{ C } \\
&= \sum_{ C \in \C_k( G )}
\sum_{ \substack{ \pi \in \Pi(C) \\
\pi'(\pi) = \{ C_1, \ldots, C_n \} } }
\sum_{ H \in \CS( K_{|\pi|} ) } (-1)^{e_H}
\prod_{i \in [n]} (-1)^{|C_i|} \Mu{C_i}.
\end{align*}
{Since
$\Ga_{\emptyset}( G )$
is the set of all clusters of $G$ with pairwise disjoint polymers,
we then rearrange the partitions according to their corresponding polymer partitions 
and have that}
\begin{align*}
\left\{ \pi \in \Pi(C) : C \in \C_k( G ) \right\}
= \left\{ \pi \in \Pi': \{ C_1, \ldots, C_n \} \in \Ga_{\emptyset}( G ): 
\sum_{i\in[n]} |C_i| \in [k] \right\},
\end{align*}
where
\begin{align*}
\Pi' := \left\{ \pi \in \Pi( \cup_{i\in[n]} C_i ): \pi'(\pi) = (C_1, \ldots, C_n) \right\}
\end{align*}
denotes the set of partitions of $\cup_{i\in[n]} C_i$
for a given set of polymers $\{ C_1, \ldots, C_n \} \in \Ga_{\emptyset}( G )$.
Hence
\begin{align}
\sum_{ C \in \C_k( G )} (-1)^{|C|} \k(C)
= \sum_{ \substack{\{ C_1, \ldots, C_n \} \in \Ga_{\emptyset}( G ) \\ 
\sum_{i\in[n]} |C_i| \in [k] } }
\sum_{\pi \in \Pi' }
\sum_{ H \in \CS( K_{|\pi|} ) } (-1)^{e_H}
\prod_{i \in [n]} (-1)^{|C_i|} \Mu{C_i}.
\label{lhs}
\end{align}
{
Note that 
we have 
\begin{align*}
\prod_{P \in \pi} (-1)^{|P|} \Mu{ P }
= \prod_{i \in [n]} (-1)^{|C_i|} \Mu{C_i}
\end{align*}
}
if and only if 
every element of the partition $\pi \in \Pi'$ is an independent set of $\GG( C_1, \ldots, C_n )$.
Then by comparing {equations} \eqref{rhs} and \eqref{lhs}, 
it suffices to show that for all $(C_1, \ldots, C_n) \in \Ga_{\emptyset}( G )$,
{
\begin{align*}
\sum_{\pi \in \Pi' }
\sum_{ G \in \CS( K_{|\pi|} ) } (-1)^{e_G} 
\prod_{P \in \pi} \I{ P \in \Ind( \GG( C_1, \ldots, C_n ) ) }
= \sum_{ H \in \CS( \GG( C_1, \ldots, C_n ) ) } (-1)^{e_H},
\end{align*}
which follows from Lemma \ref{lm-ind-part-ursell}.
}
\end{proof}

\section{Computations {for $H_3(n, p)$}}

The goal of this section is 
to compute {the terms in Theorem \ref{tm-linear}} explicitly to prove Corollary \ref{co-3linear}.

\begin{proof}[Proof of Corollary \ref{co-3linear}]

{
For $3$-uniform hypergraphs, the forbidden hypergraph 
is on four vertices with a pair of $3$-sets sharing two vertices,
we call it a {\it link}. 
Then for random indicators of links,
we construct the dependency graph $D$ following {equation} \eqref{Dind},
such that two links are adjacent if and only if they share one hyperedge.
A polymer { $C \in \C(D)$ } of size $k$ is a set of links 
$\{ F_1, \ldots, F_k \}$ 
whose induced subgraph in $D$ is connected.

{
We first enumerate all contributing non-isomorphic types of clusters,
and compute value
$\phi(\ga) (-1)^{\|\ga\|}  \allowbreak \prod_{C \in \ga} \Mu{C} / |\ga|!$ for each cluster type $\ga$.
Then we multiply each value with the size of the respective isomorphism class.
More precisely,
noting a cluster is a set of link sets, 
an isomorphism between two clusters $\ga_1, \ga_2$ is a bijection
between their vertices (the union of vertices in all links):
{$\cup_{C \in \ga_1} \cup_{F \in C}  V(F)
\to \cup_{C \in \ga_2} \cup_{F \in C} V(F)$},
which induces a bijection from the hyperedges of $\ga_1$ to the hyperedges of $\ga_2$, 
and a bijection from the polymers of $\ga_1$ to the polymers of $\ga_2$.
An automorphism of a cluster is an isomorphism to itself.
For each cluster $\ga \in \Ga(D)$, 
we consider the distinct copies of $\ga$ in the complete $r$-graph on $n$ vertices
by choosing all the vertices in 
{$\cup_{C \in \ga_1} \cup_{F \in C}  V(F)$}
from $[n]$ (ordered selections without repetition),
and every element of $\Ga(D)$ isomorphic to $\ga$ is counted once for every automorphism of $\ga$.}

Now, we compute the terms in \eqref{lprop} for $p = \o{ n^{-7/5} }$ explicitly.
}

\begin{enumerate}[leftmargin=1cm, label=(c\arabic*)]
\item Clusters $\ga$ such that $\| \ga \| = 1$.

There is only one cluster type, a single forbidden link,
namely, a hypergraph with two hyperedges intersecting in two vertices.
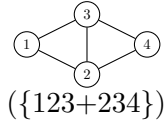
\begin{figure}[H]
\begin{center}
\begin{tikzpicture}[scale=0.8]
\tikzstyle{v}=[draw,circle,scale=0.5]
\node[v] (t21) at (2,0.5) {1};
\node[v] (t22) at (3,0) {2};
\node[v] (t23) at (3,1) {3};
\node[v] (t24) at (4,0.5) {4};
\node () at (3,-.5) { (\{123+234\}) };
\draw (t21) -- (t22) -- (t23) -- (t21);
\draw (t22) -- (t24) -- (t23); 
\end{tikzpicture}
\end{center}
\end{figure}
Thus, we have that
\begin{align*}
L^{\emptyset}_{D, 1}
= - \sum_{ C \in \C(D): |C| = 1 } \Mu{C}
= - \frac{[n]_4 p^2}{4}
= {- \frac{1}{4} n^4 p^2 + \frac{3}{2} n^3 p^2 + o(1)}.
\end{align*}

\item Clusters $\ga$ such that $\| \ga \| = 2$.

There are two cluster types: one polymer of size two, and two polymers of size one,
namely, one polymer consisting of two edge-sharing forbidden links, 
or two edge-sharing polymers with each being a single forbidden link.
Also note that for any two {(not necessarily distinct)} 
polymers $ ( C_i, C_j ) \in \C(D)^2$ such that $C_i \sim C_j$,
we have that $\GG( C_i, C_{j} ) = K_{2}$, thus $\phi(C_i, C_{j}) = -1$.

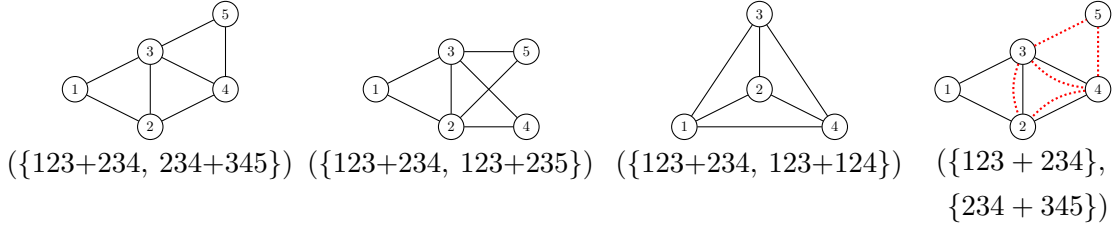
\begin{figure}[H]
\begin{center}
\begin{tikzpicture}
\tikzstyle{v}=[draw,circle,scale=0.5]
\node[v] (t31) at (4.4,0.5) {1};
\node[v] (t32) at (5.4,0) {2};
\node[v] (t33) at (5.4,1) {3};
\node[v] (t34) at (6.4,0.5) {4};
\node[v] (t35) at (6.4,1.5) {5};
\node () at (5.4,-.5) { (\{123+234, 234+345\}) };
\draw (t31) -- (t32) -- (t33) -- (t31);
\draw (t32) -- (t34) -- (t33); 
\draw (t33) -- (t35) -- (t34); 

\node[v] (t41) at (8.4,0.5) {1};
\node[v] (t42) at (9.4,0) {2};
\node[v] (t43) at (9.4,1) {3};
\node[v] (t44) at (10.4,0) {4};
\node[v] (t45) at (10.4,1) {5};
\node () at (9.4,-.5) { (\{123+234, 123+235\}) };
\draw (t41) -- (t42) -- (t43) -- (t41);
\draw (t42) -- (t44) -- (t43); 
\draw (t42) -- (t45) -- (t43); 

\node[v] (t51) at (12.5,0) {1};
\node[v] (t52) at (13.5,0.5) {2};
\node[v] (t53) at (13.5,1.5) {3};
\node[v] (t54) at (14.5,0) {4};
\node () at (13.5,-.5) { (\{123+234, 123+124\}) };
\draw (t51) -- (t52) -- (t53) -- (t51);
\draw (t54) -- (t52) (t53) -- (t54);
\draw (t54) -- (t51);

\node[v] (t31) at (16,0.5) {1};
\node[v] (t32) at (17,0) {2};
\node[v] (t33) at (17,1) {3};
\node[v] (t34) at (18,0.5) {4};
\node[v] (t35) at (18,1.5) {5};
\node () at (17,-.8) {
$\begin{aligned}
(&\{123+234\},\\
&\{234+345\})
\end{aligned}$ 
};
\draw (t31) -- (t32) -- (t33) -- (t31);
\draw (t32) -- (t34) -- (t33);
\draw[densely dotted, red, thick, bend right=20] 
(t33) edge (t32);
\draw[densely dotted, red, thick, bend right=20] 
(t33) edge (t34);
\draw[densely dotted, red, thick, bend right=20] 
(t34) edge (t32); 
\draw[densely dotted, red, thick] 
(t33) edge (t35);
\draw[densely dotted, red, thick] 
(t34) edge (t35); 

\end{tikzpicture}
\end{center}
\end{figure}
Therefore, we get
\begin{align*}
L^{\emptyset}_{D, 2}
&= \sum_{ \ga \in \Ga_{\emptyset}( D ): \|\ga\| = 2 }
\frac{ \phi(\ga) }{ |\ga|! } (-1)^{\|\ga\|}
\prod_{C \in \ga} \Mu{C} \\
&= \sum_{ C \in \C(D): |C| = 2 } (-1)^{|C|} \Mu{C}
+
\sum_{ \substack{ ( C_1, C_2 ) \in \Ga_{\emptyset}( D )  \\
|C_1| = |C_2| = 1 } }
- \frac{ 1 }{ 2 } (-1)^{ 2 } \Mu{C_1} \Mu{C_2} \\
&= \frac{[n]_5 p^3}{2} + \frac{[n]_5 p^3}{4} + \frac{[n]_4 p^3}{2} 
- \frac{ [n]_{5} p^{4} }{ 4 }
= { \frac{3}{4} n^5 p^3 + \o{1} }.
\end{align*}

\item Clusters $\ga$ such that $\| \ga \| = 3$.

We only focus on one cluster type: one polymer of size three,
namely, one polymer consisting of three edge-sharing forbidden links,
since if the cluster is formed by more then one polymer, 
then it must be extended from clusters $\ga$ such that $\| \ga \| = 2$ 
{and more than one polymer},
which are already asymptotically negligible.

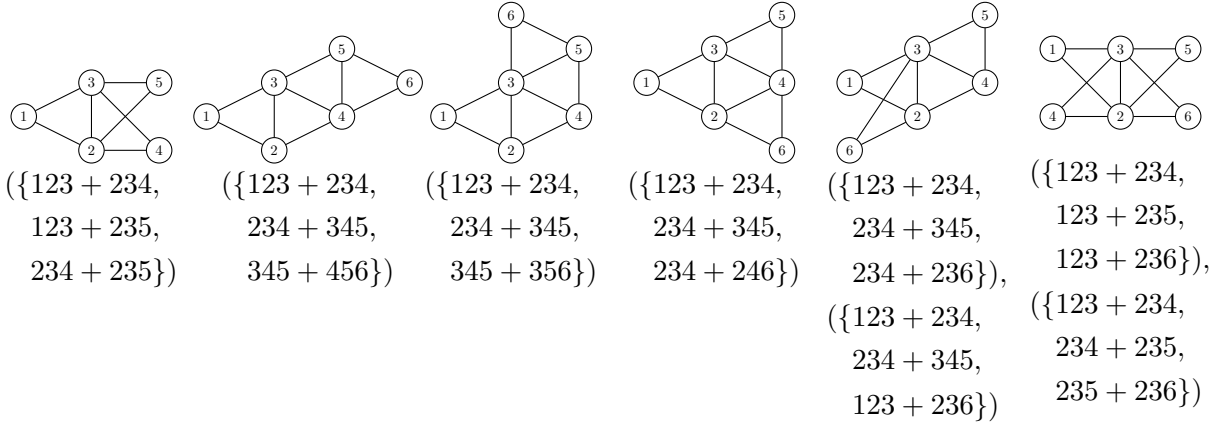
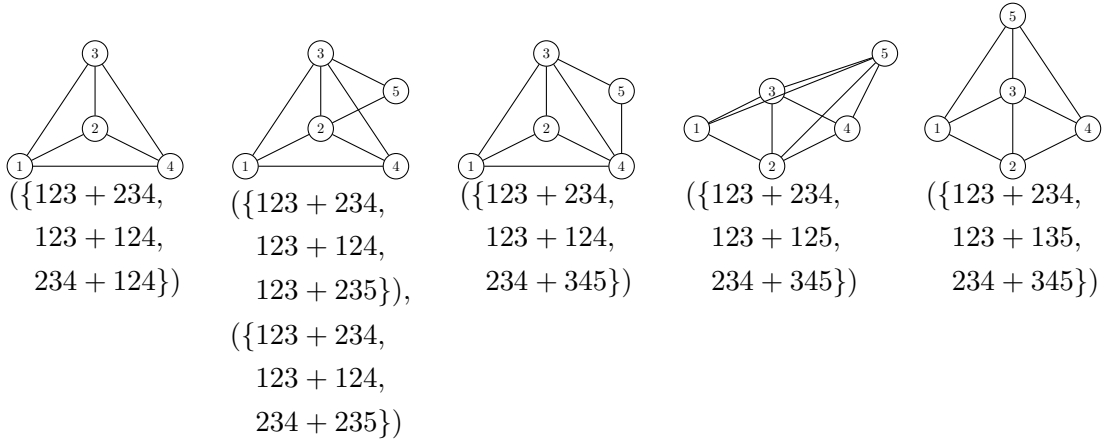
\begin{figure}[H]
\begin{center}
\begin{tikzpicture}[scale=0.9]
\tikzstyle{v}=[draw,circle,scale=0.5]
\node[v] (t41) at (2.8,0.5) {1};
\node[v] (t42) at (3.8,0) {2};
\node[v] (t43) at (3.8,1) {3};
\node[v] (t44) at (4.8,0) {4};
\node[v] (t45) at (4.8,1) {5};
\node () at (3.8,-1.2) {
$\begin{aligned}
(\{&123+234, \\&123+235, \\&234+235\})
\end{aligned}$ 
};
\draw (t41) -- (t42) -- (t43) -- (t41);
\draw (t42) -- (t44) -- (t43); 
\draw (t42) -- (t45) -- (t43); 

\node[v] (t51) at (5.5,0.5) {1};
\node[v] (t52) at (6.5,0) {2};
\node[v] (t53) at (6.5,1) {3};
\node[v] (t54) at (7.5,0.5) {4};
\node[v] (t55) at (7.5,1.5) {5};
\node[v] (t56) at (8.5,1) {6};
\node () at (7,-1.2) {
$\begin{aligned}
(\{&123+234, \\&234+345, \\&345+456\})
\end{aligned}$ 
};
\draw (t51) -- (t52) -- (t53) -- (t51);
\draw (t52) -- (t54) -- (t53); 
\draw (t53) -- (t55) -- (t54); 
\draw (t54) -- (t56) -- (t55); 

\node[v] (t61) at (9,0.5) {1};
\node[v] (t62) at (10,0) {2};
\node[v] (t63) at (10,1) {3};
\node[v] (t64) at (11,0.5) {4};
\node[v] (t65) at (11,1.5) {5};
\node[v] (t66) at (10,2) {6};
\node () at (10,-1.2) {
$\begin{aligned}
(\{&123+234, \\&234+345, \\&345+356\})
\end{aligned}$ 
};
\draw (t61) -- (t62) -- (t63) -- (t61);
\draw (t62) -- (t64) -- (t63); 
\draw (t63) -- (t65) -- (t64); 
\draw (t63) -- (t66) -- (t65); 

\node[v] (t81) at (12,1) {1};
\node[v] (t82) at (13,0.5) {2};
\node[v] (t83) at (13,1.5) {3};
\node[v] (t84) at (14,1) {4};
\node[v] (t85) at (14,2) {5};
\node[v] (t86) at (14,0) {6};
\node () at (13,-1.2) {
$\begin{aligned}
(\{&123+234, \\&234+345, \\&234+246\})
\end{aligned}$ 
};
\draw (t81) -- (t82) -- (t83) -- (t81);
\draw (t82) -- (t84) -- (t83); 
\draw (t83) -- (t85) -- (t84); 
\draw (t82) -- (t86) -- (t84); 

\node[v] (t31) at (15,1) {1};
\node[v] (t32) at (16,0.5) {2};
\node[v] (t33) at (16,1.5) {3};
\node[v] (t34) at (17,1) {4};
\node[v] (t35) at (17,2) {5};
\node[v] (t36) at (15,0) {6};
\node () at (16,-2.2) {
$\begin{aligned}
(\{&123+234, \\&234+345, \\&234+236\}), \\
(\{&123+234, \\&234+345, \\&123+236\})
\end{aligned}$ 
};
\draw (t31) -- (t32) -- (t33) -- (t31);
\draw (t32) -- (t34) -- (t33); 
\draw (t33) -- (t35) -- (t34); 
\draw (t32) -- (t36) -- (t33); 

\node[v] (t41) at (18,1.5) {1};
\node[v] (t42) at (19,0.5) {2};
\node[v] (t43) at (19,1.5) {3};
\node[v] (t44) at (18,.5) {4};
\node[v] (t45) at (20,1.5) {5};
\node[v] (t46) at (20,.5) {6};
\node () at (19,-2) {
$\begin{aligned}
(\{&123+234, \\&123+235, \\&123+236\}), \\
(\{&123+234, \\&234+235, \\&235+236\})
\end{aligned}$ 
};
\draw (t41) -- (t42) -- (t43) -- (t41);
\draw (t42) -- (t44) -- (t43); 
\draw (t42) -- (t45) -- (t43); 
\draw (t42) -- (t46) -- (t43); 

\end{tikzpicture}
\end{center}
\smallskip
\begin{center}
\begin{tikzpicture}
\tikzstyle{v}=[draw,circle,scale=0.5]

\node[v] (t51) at (-2,0) {1};
\node[v] (t52) at (-1,0.5) {2};
\node[v] (t53) at (-1,1.5) {3};
\node[v] (t54) at (0,0) {4};
\node () at (-1,-1) {
$\begin{aligned}
(\{&123+234, \\&123+124, \\&234+124\})
\end{aligned}$ 
};
\draw (t51) -- (t52) -- (t53) -- (t51);
\draw (t54) -- (t52) (t53) -- (t54);
\draw (t54) -- (t51);

\node[v] (t51) at (1,0) {1};
\node[v] (t52) at (2,0.5) {2};
\node[v] (t53) at (2,1.5) {3};
\node[v] (t54) at (3,0) {4};
\node[v] (t55) at (3,1) {5};
\node () at (2,-2) {
$\begin{aligned}
(\{&123+234, \\&123+124, \\&123+235\}), \\
(\{&123+234, \\&123+124, \\&234+235\}) \\
\end{aligned}$ 
};
\draw (t51) -- (t52) -- (t53) -- (t51);
\draw (t54) -- (t52) (t53) -- (t54);
\draw (t54) -- (t51);
\draw (t52) -- (t55) -- (t53);

\node[v] (t51) at (4,0) {1};
\node[v] (t52) at (5,0.5) {2};
\node[v] (t53) at (5,1.5) {3};
\node[v] (t54) at (6,0) {4};
\node[v] (t55) at (6,1) {5};
\node () at (5,-1) {
$\begin{aligned}
(\{&123+234, \\&123+124, \\&234+345\})
\end{aligned}$ 
};
\draw (t51) -- (t52) -- (t53) -- (t51);
\draw (t54) -- (t52) (t53) -- (t54);
\draw (t54) -- (t51);
\draw (t53) -- (t55) -- (t54);

\node[v] (t31) at (7,0.5) {1};
\node[v] (t32) at (8,0) {2};
\node[v] (t33) at (8,1) {3};
\node[v] (t34) at (9,0.5) {4};
\node[v] (t35) at (9.5,1.5) {5};
\node () at (8,-1) {
$\begin{aligned}
(\{&123+234, \\&123+125, \\&234+345\})
\end{aligned}$ 
};
\draw (t31) -- (t32) -- (t33) -- (t31);
\draw (t32) -- (t34) -- (t33); 
\draw (t33) -- (t35) -- (t34); 
\draw (t32) -- (t35) -- (t31); 

\node[v] (t31) at (10.2,0.5) {1};
\node[v] (t32) at (11.2,0) {2};
\node[v] (t33) at (11.2,1) {3};
\node[v] (t34) at (12.2,0.5) {4};
\node[v] (t35) at (11.2,2) {5};
\node () at (11.2,-1) {
$\begin{aligned}
(\{&123+234, \\&123+135, \\&234+345\})
\end{aligned}$ 
};
\draw (t31) -- (t32) -- (t33) -- (t31);
\draw (t32) -- (t34) -- (t33); 
\draw (t31) -- (t35) -- (t33);
\draw (t35) -- (t34); 
\end{tikzpicture}
\end{center}
\end{figure}

Hence, we have that
\begin{align*}
L^{\emptyset}_{D, 3}
&= \sum_{ \ga \in \Ga_{\emptyset}( D ): \|\ga\| = 3 }
\frac{ \phi(\ga) }{ |\ga|! } (-1)^{\|\ga\|}
\prod_{C \in \ga} \Mu{C} \\
&= \sum_{ C \in \C( D ): |C| = 3 } (-1)^{|C|} \Mu{C}
+ \O{ n^4 p^3 } + \O{ n^5 p^4 } \\
&= - \frac{[n]_5 p^3}{2 \times 3!} 
- \frac{[n]_6 p^4 }{2}  
- \frac{[n]_6 p^4 }{2}
- \frac{[n]_6 p^4 }{3!}
- \frac{[n]_6 p^4 }{2}
- [n]_6 p^4
- \frac{[n]_6 p^4 }{2 \times 3!} 
- \frac{[n]_6 p^4 }{{2 \times 2}} 
+ o(1),
\end{align*}
where the last row of the types of polymers are of contribution $\O{ n^4 p^3 } = o(1)$.

{

\item Clusters $\ga$ such that $\| \ga \| = 4$. \\
As before, we only focus on one cluster type: one polymer of size four,
since clusters with more than one polymer contribute negligibly.
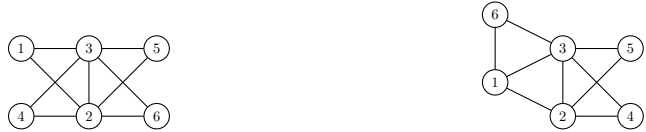
\begin{figure}[H]
\begin{center}
\begin{tikzpicture}[scale=0.9]
\tikzstyle{v}=[draw,circle,scale=0.5]

\node[v] (t46) at (6,2) {6};
\node[v] (t41) at (6,1) {1};
\node[v] (t42) at (7,0.5) {2};
\node[v] (t43) at (7,1.5) {3};
\node[v] (t44) at (8,.5) {4};
\node[v] (t45) at (8,1.5) {5};

\node () at (7,-.5) {
$\begin{aligned}
(\{ 123+234, 123+235, \\
234+235, 123+136\})
\end{aligned}$ 
};
\draw (t41) -- (t42) -- (t43) -- (t41);
\draw (t42) -- (t44) -- (t43); 
\draw (t42) -- (t45) -- (t43);
\draw (t41) -- (t46) -- (t43); 

\node[v] (t41) at (-1,1.5) {1};
\node[v] (t42) at (0,0.5) {2};
\node[v] (t43) at (0,1.5) {3};
\node[v] (t44) at (-1,.5) {4};
\node[v] (t45) at (1,1.5) {5};
\node[v] (t46) at (1,.5) {6};
\node () at (0,-.5) {
$\begin{aligned}
&(\{123+234, 123+235, 123+236, 234+235\}), \\
&(\{123+234, 234+235, 235+236, 123+236\})
\end{aligned}$};
\draw (t41) -- (t42) -- (t43) -- (t41);
\draw (t42) -- (t44) -- (t43); 
\draw (t42) -- (t45) -- (t43); 
\draw (t42) -- (t46) -- (t43); 
\end{tikzpicture}
\end{center}
\end{figure}
We then have
\begin{align*}
L^{\emptyset}_{D, 4}
&= \sum_{ \ga \in \Ga_{\emptyset}( D ): \|\ga\| = 4 }
\frac{ \phi(\ga) }{ |\ga|! } (-1)^{\|\ga\|}
\prod_{C \in \ga} \Mu{C} 
= \frac{[n]_6 p^4 }{2 \times 2} 
+ \frac{[n]_6 p^4 }{2 \times 8}
+ \frac{[n]_6 p^4 }{2} + o(1).
\end{align*}

\item Clusters $\ga$ such that $\| \ga \| \in \{ 5, 6 \}$.
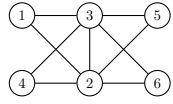
\begin{figure}[H]
\begin{center}
\begin{tikzpicture}[scale=0.9]
\tikzstyle{v}=[draw,circle,scale=0.5]
\node[v] (t41) at (-1,1.5) {1};
\node[v] (t42) at (0,0.5) {2};
\node[v] (t43) at (0,1.5) {3};
\node[v] (t44) at (-1,.5) {4};
\node[v] (t45) at (1,1.5) {5};
\node[v] (t46) at (1,.5) {6};
\node () at (0,-.5) {
$\begin{aligned}
(\{123+234&, 123+235, \\
123+236, 234&+235, 234+236 \}) \\
\end{aligned}$};
\draw (t41) -- (t42) -- (t43) -- (t41);
\draw (t42) -- (t44) -- (t43); 
\draw (t42) -- (t45) -- (t43); 
\draw (t42) -- (t46) -- (t43); 

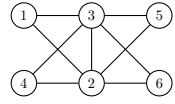
\node[v] (t41) at (6,1.5) {1};
\node[v] (t42) at (7,0.5) {2};
\node[v] (t43) at (7,1.5) {3};
\node[v] (t44) at (6,.5) {4};
\node[v] (t45) at (8,1.5) {5};
\node[v] (t46) at (8,.5) {6};
\node () at (7,-.5) {
$\begin{aligned}
(\{&123+234, 123+235, 123+236, \\
&234+235, 234+236, 235+236\})
\end{aligned}
$};
\draw (t41) -- (t42) -- (t43) -- (t41);
\draw (t42) -- (t44) -- (t43); 
\draw (t42) -- (t45) -- (t43); 
\draw (t42) -- (t46) -- (t43); 
\end{tikzpicture}
\end{center}
\end{figure}
Then we have
\begin{align*}
L^{\emptyset}_{D, 5} + L^{\emptyset}_{D, 6}
&= \sum_{ \ga \in \Ga_{\emptyset}( D ): \|\ga\| \in \{5, 6\} }
\frac{ \phi(\ga) }{ |\ga|! } (-1)^{\|\ga\|}
\prod_{C \in \ga} \Mu{C} 
= - \frac{[n]_6 p^4 }{2 \times 2 \times 2}
+ \frac{[n]_6 p^4 }{2 \times 4!} 
+ o(1).
\end{align*}
}
\end{enumerate}

{
Since there is a finite number of types of polymers with size seven,
we thus have $\Delta_{7}(D) = \o{1}$.
Hence, we ignore the remaining terms by equation \eqref{exp-del}.
}
Adding up the contributing terms for $p = { \o{ n^{-7/5} } }$
gives the asymptotic probability of $H_3(n, p)$ being linear in Corollary \ref{co-3linear}.
\end{proof}

\section{Concluding remarks}


We have shown that the truncation of the cluster expansion series gives 
the asymptotic linearity of binomial random hypergraphs.
The analysis of the truncation utilised results in \cite{mousset2020probability},
which exploit the correlation among random variables 
and rely heavily on FKG inequality.
It would be interesting to investigate whether 
this is necessary for the truncation.
Alternative ways of handling truncation that are commonly used include 
establishing the absolute convergence of the series 
via the Koteck\`y-Preiss criterion \cite{kotecky1986cluster},
for example, see \cite{helmuth2020algorithmic,jenssen2020algorithms}.

It also would be interesting to investigate whether
the cluster expansion series also gives the linearity of random hypergraphs with given number of edges.
However, in that case, all graph-dependent indicators are dependent,
and the only valid dependency graph for them is {the} complete graph.
In this case, we may need to modify the method by incorporating the notion of weak dependence,
see, for example, \cite{isaev2021extremal}.

\section*{Acknowledgements}
{The author sincerely thanks Brendan McKay for significant help,
including checking {answers} via numerical simulations, 
providing his computational results,
pointing out calculation errors, and many email discussions.}
The author would like to thank his supervisor Nick Wormald 
for introducing the problem,
many critical discussions, and comments.
The author also thanks Will Perkins for critical discussions on the cluster expansion.
The author is sincerely grateful
to an anonymous referee for
{carefully reading the manuscript,}
pointing out errors,
and providing invaluable comments and suggestions, 
which lead to a substantial improvement of the presentation {and quality}.

\bibliographystyle{plain}
\bibliography{ref}

\end{document}